\documentclass[12pt]{amsart}

\usepackage{graphicx}
\usepackage{amsthm}
\swapnumbers

\newtheorem{Thm}{Theorem}
\newtheorem{MThm}[Thm]{Main Theorem}
\newtheorem{Def}[Thm]{Definition}
\newtheorem{Coro}[Thm]{Corollary}
\newtheorem{Lem}[Thm]{Lemma}
\newtheorem{Sublem}[Thm]{Sublemma}
\newtheorem{Rem}[Thm]{Remark}
\begin{document}

\title{Knots with infinitely many incompressible Seifert surfaces}
\author{Robin T. Wilson}

\begin{abstract}
We show that a knot in $S^3$ with an infinite number of incompressible Seifert surfaces contains a closed incompressible surface in its complement.
\end{abstract}

\maketitle

\section{Introduction}
\label{defsect}

In his 1978 thesis, Parris \cite{Pa} showed that there exist pretzel knots with an infinite number of incompressible Seifert surfaces.  This result was later used by Casson and Gordon \cite{CG} to construct a 3-manifold with an infinite number of non-isotopic strongly irreducible Heegaard splittings of arbitrarily high genus, contradicting a conjecture of Waldhausen \cite{Wa}.  While studying these examples, Kobayashi \cite{Ko2} observed that the entire family of incompressible Seifert surfaces used in the Casson-Gordon example is carried by one incompressible branched surface.  In a similar spirit, Moriah, Sedgwick, and Schleimer \cite{MSS} proved that this same family of surfaces can be constructed by taking a Haken sum of one of the Seifert surfaces in the family with multiple copies of a closed incompressible surface of genus 2.  Tao Li \cite{L1} observed, using branched surfaces, that this incompressible surface of genus 2 can be considered to be a limit of a sequence of strongly irreducible Heegaard surfaces.  

The goal of this paper is to prove a related result about Seifert surfaces in knot complements.  We show that every knot with an infinite number of distinct incompressible Seifert surfaces contains a closed incompressible surface in its complement.  Yoav Rieck and Eric Sedgwick have pointed out that they have included a brief sketch of a proof of this in \cite{RS}.  We do this by showing that, in a knot complement, all incompressible Seifert surfaces can be constructed from a finite list of closed incompressible surfaces and incompressible Seifert surfaces for the knot as a linear combination using the Haken sum.  Our main theorem is the following:

\begin{MThm}
\label{maintheorem}
Let $K$ be a non-trivial knot in $S^3$ and let $M=S^3-N(K)$.  Then there are a finite number of incompressible Seifert surfaces $\{S_1,...,S_n\}$ and a finite number of closed incompressible surfaces $\{Q_1,...,Q_m\}$ in $M$, none of which are boundary parallel, such that any incompressible Seifert surface $S$ is isotopic to a Haken sum $S=S_i+a_1Q_1+\cdots+a_mQ_m$, where the $a_1,...,a_m$ are non-negative integers. 
\end{MThm}

\begin{Coro}  Every knot with an infinite number of distinct incompressible Seifert surfaces contains a closed incompressible surface in its complement.  
\end{Coro}
\begin{proof}
Suppose that $K$ is a knot with an infinite number of distinct incompressible Seifert surfaces $\{G_i\}_{i=1}^{\infty}$.  Then by Theorem ~\ref{maintheorem} there exists a finite set of incompressible Seifert surfaces $\{S_1,...,S_n\}$ and a finite set of closed incompressible surfaces $\{Q_1,...,Q_m\}$ (possibly empty), such that for each $i\geq1$, $G_i$ can be isotoped so that $G_i=S_j+a_1Q_1+\cdots+a_mQ_m$.  Since we have an infinite collection on the left hand side of this equation and only a finite number of $S_j$'s and $Q_k$'s, the only possibility is that the set $\{Q_1,...,Q_m\}$ is non-empty.  Thus for some $l \in \{1,...,m\}$, $Q_l$ is a closed incompressible surface in the complement of $K$.
\end{proof}

 The paper is organized as follows.  In Section 2 we review normal surface theory from the viewpoint of triangulations.  Given a triangulation of a 3-manifold, one obtains a system of linear equations and inequalities called the normal surface equations.  Each solution to these equations corresponds to a normal surface in the 3-manifold.   In Section 3 we introduce a modification to the normal surface equations that restricts the class of surfaces that arise as solutions to the normal surface equations in a knot complement to Seifert surfaces and closed surfaces only.  In Section 4 we use the modified normal surface theory to prove our main theorem using a result of \cite{JO}.  The paper concludes in Section 5 with an explanation of the motivating example of Parris \cite{Pa} of a pretzel knot with an infinite family of incompressible Seifert surfaces.  We look at how this family can be expressed as the Haken sum of one incompressible Seifert surface with multiple copies of a closed incompressible surface. \\

I would like to thank Joel Hass, Abigail Thompson, and Hyam Rubinstein for their helpful conversations.  I also would like to thank Jesse Johnson for his many valuable comments.  This research was done while under the support of the National Academy of Sciences.

\section{Review of Normal Surface Theory}
\label{normsurfsect}

Throughout the paper we will assume that $M$ is an orientable 3-manifold, and that $S$ is a properly embedded, compact, orientable surface.  If $K$ is a knot in $S^3$ then we will assume that $K$ is non-trivial and we will let $N(K)$ denote a regular neighborhood of $K$ in $S^3$.

\begin{Def}\textup{
Let $S$ be a triangulated surface and let $c$ be a curve on  $S$.  Assume that $c$ is transverse to the 1-skeleton of the triangulation.  A curve $c$ in $S$ is called \textit{normal} if the intersection of $c$ with any triangle of the triangulation contains no closed curves and no arcs with both endpoints on the same edge.
}\end{Def}

\begin{Def}\textup{
Let $M$ be a triangulated 3-manifold.  A \textit{normal triangle} in a tetrahedron of the triangulation is an embedded disk that meets three edges and three faces of the tetrahedron.  A  \textit{normal quadrilateral} is an embedded disk in a tetrahedron that meets four edges and four faces of the tetrahedron.  Normal triangles and quadrilaterals are called  \textit{normal disks}.
}\end{Def}  

\begin{figure}
  \begin{center}
  \includegraphics[width=2.5in]{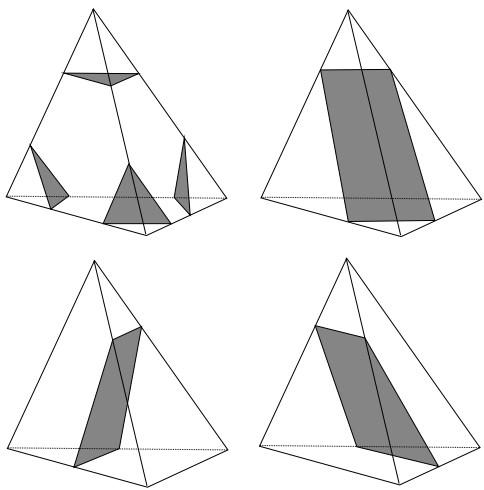}
  \caption{Normal triangles and quadrilaterals}
  \label{normaldisks}
  \end{center}
\end{figure}

\begin{Def}\textup{
Let $M$ be a triangulated 3-manifold.  An embedded surface $S\subset M$ is a \textit{normal surface} if it meets each tetrahedron in a disjoint collection of normal disks. 
}\end{Def}

\begin{Def}\textup{
A surface $S$ embedded in a 3-manifold is \textit{compressible} if either $S$ is a 2-sphere bounding a 3-ball in $M$, or there is an essential, simple, closed curve $c$ in $S$ which bounds a disk $D$ embedded in $M$ such that $int(D)  \cap S= \emptyset$.  The disk $D$ is called a \textit{compressing disk} for $S$.  A surface that is not compressible is called \textit{incompressible}.
}\end{Def}

\begin{Def}\textup{
The \textit{complexity} of a normal surface $S$ in a triangulated 3-manifold M, denoted $\gamma(S)$, is the number of times $S$ intersects the 1-skeleton of the triangulation.  A normal surface of minimal complexity in its isotopy class is called \textit{minimal}.
}\end{Def}

The following theorems are due to Haken \cite{H}:

\begin{Thm}
\label{incompthm}
Let $M$ be an irreducible 3-manifold containing an incompressible surface $S$.  Then for any triangulation of $M$ there is an isotopy of $S$ that takes $S$ to a normal surface in $M$.  
\end{Thm}

\begin{Thm}
\label{normalcurvethm}
Any non-trivial simple closed curve $c$ in a triangulated surface $S$ can be isotoped to be normal with respect to the given triangulation.  
\end{Thm}

Let $M$ be a 3-manifold with a given triangulation.  For each tetrahedron there are four types of normal triangles and three types of normal quadrilaterals as shown in Figure \ref{normaldisks}.  To describe a normal surface in $M$ up to isotopy we only need to specify how many of each type of normal disks occurs in each tetrahedron.  If the triangulated 3-manifold $M$ has $t$ tetrahedra then a normal surface is completely described by a vector of $7t$ non-negative integers.  

Conversely, for a non-negative integer vector in $Z^{7t}$ to give rise to an embedded normal surface two conditions must be met:

\begin{itemize}
\item Each tetrahedron must contain at most one type of quadrilateral.  This is called the \emph{quadrilateral restriction}.
\item The normal disks must match up across each face shared by two adjacent tetrahedra.  

\end{itemize}
The first condition places a restriction on the number of non-zero coordinates that appear in a vector representing a normal surface.  This is because two different types of normal quadrilaterals in a tetrahedron must intersect and we want to avoid this to ensure that our normal surface will be embedded.  The second condition becomes necessary when one considers two adjacent tetrahedra sharing a face.  The number of normal disks in one tetrahedron yielding arcs in the common face must equal the number of normal disks meeting the same face in the adjacent tetrahedron.  This gives a sequence of linear equations of the form $x_i+x_j=x_k+x_l$.  These linear equations are called the \textit{matching equations}.  Note that if the 3-manifold has boundary then there are no matching equations associated to the boundary faces.  We add one more restriction, that all solutions must be non-negative and integral.

This system of equations is called the \textit{system of normal surface equations}. If a vector is a solution to the normal surface equations then the given vector corresponds to a unique normal surface in $M$ up to isotopy. 

\begin{figure}
  \begin{center}

  \includegraphics[width=2.5in]{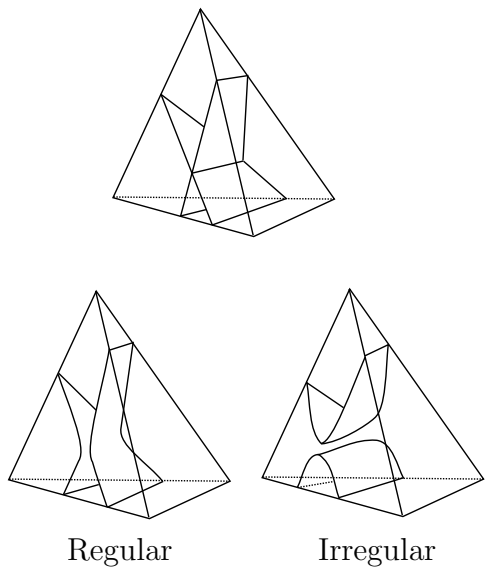}
  \put(-158,-15){Regular}
  \put(-63,-15){Irregular}
  \caption{The Haken sum $F+G$ }
  \label{hakensumdef}
 \end{center}
\end{figure}

Two normal surfaces $F$ and $G$ in $M$ that intersect are called \textit{compatible} with respect to a triangulation $T_M$ if in each 3-simplex, the quadrilaterals of $F$ and $G$, if any, are of the same type.  This condition ensures that there is always a way to cut and paste two normal disks intersecting along an arc to obtain a pair of disjoint normal disks.  If $F$ and $G$ are compatible, then for each component $\alpha$ of $F\cap G$,  $N(\alpha)$, a small open neighborhood of $\alpha$ in $M$,  is a solid torus (or one-handle) with boundary torus (or annulus) $T$.  The set $F\cup G$ separates this torus (annulus) $T$ into four annuli (rectangles):  $A_1, A_2, A_3, A_4$, labeled cyclically.  Let $A_+=A_1 \cup \ A_3$ and $A_-=A_2 \cup \ A_4$.  For each $\alpha$ in $F \cap G$ we remove $(F \cup G) \cap N(\alpha)$, and replace it with a choice of $A_+$ or $A_-$.  Of the two pairs $A_+$ and $A_-$, one choice gives us a normal surface and the other does not.  We will call the pair that gives us a normal surface a \textit{regular pair} of annuli (rectangles) and the other an \textit{irregular pair}.  For each arc $\alpha$ of $F\cap G$ we remove the corresponding pieces of $F \cup G$ and replace them with a regular pair of annuli (rectangles).  We denote the resulting surface, called the \textit{Haken sum} by $F+G$.   See Figure \ref{hakensumdef}.  The Haken sum $F+G$ is again a normal surface, and up to isotopy is independent of the choices involved in its construction.

 Here are some useful facts about Haken sums:
 \begin{itemize}
 \item If $F=F_1+F_2$ and $F_1$, $F_2$ correspond to the non-negative integer vectors $(x_{1},...,x_{n})$, $(y_{1},...,y_{n})$ respectively, then $F$ corresponds to the vector $(x_{1}+y_{1},...,x_{n}+y_{n})$.
 \item If $F=F_1+F_2$, then $\chi (F)=\chi (F_1)+\chi(F_2)$, where $\chi$ is the Euler characteristic.
 \item If $F=F_1+F_2$, then $\gamma(F)=\gamma(F_1)+\gamma(F_2)$, where $\gamma$ is the complexity.  
 \end{itemize}
 
\begin{Def}\textup{
A solution $\vec x=(x_{1},...x_{n})$ of the system of normal surface equations is a \textit{fundamental solution} if it cannot be written as the sum of two other non-zero solutions.  That is, $\vec x \ne \vec y + \vec z$, where $\vec y, \vec z$ are any solutions to the normal surface equations. 
}\end{Def}

\begin{Def}\textup{
If $\vec x$ is a fundamental solution then the corresponding normal surface is called a \textit{fundamental surface} in $M$.  
}\end{Def}

Haken showed, using the following theorem, that there is a finite set of fundamental solutions to the system of normal surface equations, so that any solution is an integer linear combination of this finite set of fundamental solutions.  A system of linear equations with only integer solutions allowed is called a \textit{Diophantine} system of equations.  Thus the system of normal surface equations os a Diophantine system.

\begin{Thm}[\cite{H}]
\label{diophantinethm}
Let $Ax=0$ be a Diophantine system of equations.  Then there is a finite set of non-negative integral solutions which generates, using non-negative integer linear combinations, the full set of all non-negative solutions to the Diophantine system of equations.
\end{Thm}

Hence any normal surface in $M$ can be expressed as a non-negative integer linear combination of the finite set of fundamental normal surfaces for $M$.

\section{Modifying the normal surface equations}
\label{joelsect}

To use normal surfaces to study incompressible Seifert surfaces in a knot complement, we would like to restrict our attention to finding only solutions of the normal surface equations that correspond to surfaces that either have a fixed preferred longitude as boundary, or are closed.  In this section we use Theorem \ref{diophantinethm} to show that we can find a finite set of fundamental solutions to suitably restricted normal surface equations so that every normal surface that is closed or that has the fixed preferred longitude as boundary is expressible as a non-negative integer linear combination of this finite set. 

\begin{Lem}
\label{curvelem}
Let $M$ be a 3-manifold with $\partial M\neq \emptyset$, and let $\alpha$ be a fixed curve in $\partial M$.  Then there is a triangulation of $M$ such that $\alpha$ meets each 2-simplex of the triangulation of $\partial M$ in at most one normal arc.  
\end{Lem}

\begin{proof}  Given $M$ and $\alpha$ as above, let $T_M$ be a triangulation of the 3-manifold $M$.  Now consider the triangulation $T_{\partial M}$ obtained by restricting $T_M$ to $\partial M$.  By Theorem \ref{normalcurvethm},  we can isotope $\alpha$ so that it is normal with respect to $T_{\partial M}$.  Each triangle $\sigma$ in $T_M$ is isotopic to a planar equilateral triangle $\sigma'$ and the arcs of $\alpha \cap \sigma$ can be mapped to straight line segments.  We describe a barycentric subdivision for the planar equilateral triangles $\sigma'$ which we will pull back to the original triangle $\sigma$ of the triangulation $T_M$.  

Consider a fixed  triangle $\sigma$ isotopic to the planar $\sigma'$.  Then given an edge $e$ of $\sigma'$, the image of the curve $\alpha$ divides $e$ into regions.  A region of $e$ is an \textit{inner edge} if it doesn't contain a vertex of $\sigma$.  The curve $\alpha$ divides $\sigma$ into regions whose boundary consists of arcs of $\alpha$ and subarcs of edges of $\sigma$.  An \textit{inner region} of $\sigma$ is a region that has 2 or more arcs of  the image of $\alpha$ in its boundary.  (An inner region of $\sigma$ may contain a vertex).

We will slightly modify the usual barycentric subdivison process of the planar equilateral triangle $\sigma'$ by being careful about where we place the barycenters.  For each edge of $T_M$ that meets the image of $\alpha$  in more than one point, place a barycenter in an inner edge.  For each edge of $T_M$ that has no inner edges use the `true' barycenter.  For each face of $T_M$ meeting the image of $\alpha$ in more than one arc place a barycenter in an inner region.  If a face has no inner regions use the true barycenter.  

We now pull back the subdivided triangles $\sigma'$ to the original triangle $\sigma$ in $T_M$.  We will call the resulting triangulation the 1st barycentric subdivision of $T_M$, written ${T_M}^{(1)}$.  Subdividing $T_M$ in this fashion reduces the maximum number of inner regions and inner edges in each triangle of $T_{\partial M}^{(1)}$.  Thus this process reduces the maximum number of components of $\alpha \cap \sigma$ over all triangles $\sigma$ in $T_{\partial M}^{(1)}$.  We can repeatedly take barycentric subdivisions and stop when we find a triangulation of $M$, say  ${T_M}^{(N)}$ such that $\alpha$ meets each face of  $T_{\partial M}^{(N)}$ in at most one normal arc, where $N$ is some non-negative integer.  
\end{proof}

\begin{Rem}  It follows that if N is the number  of subdivisions necessary before reaching a triangulation ${T_M}^{(N)}$ such that $\alpha$ meets each face of  $T_{\partial M}^{(N)}$ in at most one normal arc, then N is less than the maximal number of arcs in $\alpha \cap \sigma$, where the maximum is taken over all 2-simplices $\sigma$ of $T_{\partial M}^{(N)}$.
\end{Rem}

\begin{Lem} 
\label{fixedarclem}  
Let $K$ be a knot in $S^3$, let $M=S^3-N(K)$ and let $S$ be an incompressible Seifert surface in $M$.  Suppose that $T_M$ is  a triangulation of $M$ such that $\partial S$ meets each 2-simplex of $T_{\partial M}$ in at most one normal arc. Then $S$ can  be made normal with respect to $T_M$ by an isotopy fixing $\partial S$. 
\end{Lem}

 \begin{proof}
This proof closely follows the proof of Theorem~\ref{incompthm} as presented in \cite{JR}.  Let $S$ be an incompressible Seifert surface in $M$.  Let $T_M$ be a triangulation of $M$ such that $\partial M$ meets each 2-simplex of $T_{\partial M}$ in at most one normal arc.  The first step in the process of normalizing $S$ is to isotope $S$ so that it is in general position with respect to the 2-skeleton of $T_M$, and so that $S$ meets the 1-skeleton of $T_M$ in the least number of components over all surfaces in its isotopy class $rel(\partial S)$ satisfying the previous condition.  
 
 The second step is to remove any simple closed curve components of the intersection of $S$ with the 2-skeleton of $T_M$.  If the intersection of $S$ with a triangle $\sigma$ of $T_M$ contains a simple closed curve $\alpha$ we may assume that the curve $\alpha$ is ``innermost" in $\sigma$.  Then $\alpha$ bounds a disk $D$ in $\sigma$ with $int(D) \cap S = \emptyset$.  Surgery on $S$ along the disk $D$ gives two components.  One is a 2-sphere bounding a 3-ball in $M$, and  the other is a surface isotopic to $S$ which we will continue to call $S$.  This step does not create any new component of intersection with the 1-skeleton.  Since surgery along $\alpha$ takes place in the interior of $M$ this step does not disturb the boundary of $S$.
 
 The third step is to show that $S$ meets each triangle of $T_M$ in normal arcs.  If some triangle $\sigma$ in $T_M$ meets $S$ in an arc $\alpha$ such that $\alpha$ has both endpoints on the same edge of $e$ of $\sigma$, then we can assume that the arc $\alpha$ is an ``outermost'' arc with both endpoints on the same edge $e$.  The edge $e$ cannot be in $\partial M$ because this would contradict the boundary condition.  So the edge $e$ must be in the interior of $M$. Then $S$ can be isotoped to remove the arc of intersection $\alpha$, reducing the number of times that $S$ meets the 1-skeleton of $T_M$ contradicting the minimality of $\gamma(S)$.   Again, after this step $\partial S$ remains fixed.  
 
 In the fourth step we show that the intersection of $S$ with each tetrahedron of $T_M$ is a disk.  Suppose $\Delta$ is a tetrahedron of $T_M$ such that some component of $S \cap \Delta$ is not a disk.  Then there is a curve $\alpha$ of $S\cap \partial \Delta$ so that $\alpha$ bounds a disk $D$ in $\Delta$, $\alpha$ does not bound a disk in $S \cap \Delta$, and $int(D)\cap S = \emptyset$.  As in step 2, surgery of $S$ along this disk $D$ gives two new surfaces, with one component a 2-sphere bounding a 3-ball, and the other a surface isotopic to $S$ but with fewer components of intersection with the 1-skeleton.  This contradicts the minimality assumption, and so we can conclude that each component of $S \cap \Delta$ is a disk.  Again, since surgery along $\alpha$ takes place in the interior of $M$, $\partial S$ remains fixed.
 
The fifth and last step is to show that each component of intersection of $S$ with the boundary of a tetrahedron is isotopic to a simple closed curve that meets each triangle on the boundary of the tetrahedron in at most one normal arc.  Again, we must consider how the process of achieving this affects the boundary condition. Let $\Delta$ be a tetrahedron of $T_M$ and suppose that there exists a component  $J$ of $S\cap \partial \Delta$ that is a simple closed curve intersecting a face $\sigma$ of the tetrahedron in more than one arc.  We may assume that $J$ is ``innermost'' in the sense that $J$ meets a face $\sigma$ of $\Delta$ in two arcs $a$ and $b$ which have two endpoints in a common edge $e$ of $\sigma$, and no other components of $S \cap e$ are contained between the endpoints $a$ and $b$.  

Above we have shown that $J$ bounds a disk $D$ in $\Delta$.  The disk $D$ co-bounds with $\sigma$ a disk $E$ that gives a boundary compression of $D$. If $e$ is in $\partial M$ then we obtain a contradiction to the boundary condition in the hypothesis. Therefore the face $\sigma$ must be in the interior of $M$,  and so the disk $E$ bounded describes an isotopy of $S$ which doesn't change the number of times that that $S$ meets the 1-skeleton, and it introduces a new component of intersection of $S$ with $\sigma$ that is an arc with both endpoints on the same edge $e$.  However this possibility has been ruled out in the second step above.  Therefore the intersection of $S$ with any tetrahedron of $T_M$ is a disk whose boundary is a curve on the boundary of the tetrahedron which meets each face of the tetrahedron in at most one arc, which must span between distinct edges of the faces.  We can conclude that $S$ has been isotoped so that it meets each tetrahedron of $T_M$ in normal disks, thus $S$ has been isotoped to be a normal surface while keeping $\partial S$ fixed.   
\end{proof}

\begin{figure}
  \begin{center}
  \includegraphics[width=2.7in]{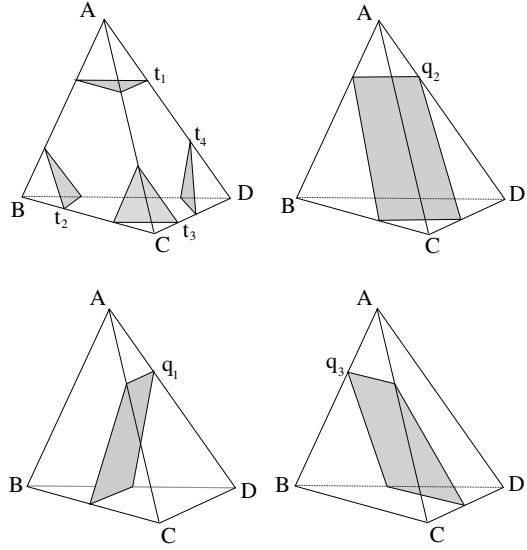}
  \caption{Notation for the normal surfaces in a tetrahedra}
  \label{tetrahedranotation}
  \end{center}
\end{figure}

In order to impose additional restrictions on the system of normal surface equations, we need some further notation.  Let $\{A,B,C,D\}$ be the four vertices determining a tetrahedron.  Assign variables $t_1,t_2,t_3$, and $t_4$ to the four normal triangles in a tetrahedron, and variables $q_1,q_2$ and $q_3$ to the three normal quadrilaterals as shown in Figure \ref{tetrahedranotation}.  Then the normal triangles are determined by specifying the corresponding vertex, and the quadrilaterals by specifying an edge not intersected by the disk.\\

\begin{itemize}
\item   $t_1, t_2, t_3,t_4 \longrightarrow A,B,C,D$
\item  $q_1,q_2,q_3 \longrightarrow AB, AC, AD$
\end{itemize}.

Let $K$ be a knot in $S^3$, and let $M=S^3-N(K)$ with triangulation $T_M$.  Let $S$ be a normal surface corresponding to a solution of the normal surface equations which has $\partial S= \alpha$, where $\alpha$ intersects each 2-simplex on $\partial M$ in at most one normal arc.  We will describe how to modify the normal surface equations so that any solution to this modified system is closed or has boundary consisting of copies of $\alpha$.  Suppose $\{A,B,C,D\}$ determines a tetrahedron of $T_M$ that has one face on $\partial M$, say the face is given by $BCD$.  Let $\sigma$ represent the face $BCD$.

For the face $\sigma$ on $\partial M$ we have two possibilities: \\

\begin{figure}
  \begin{center}
  \includegraphics[width=2.6in]{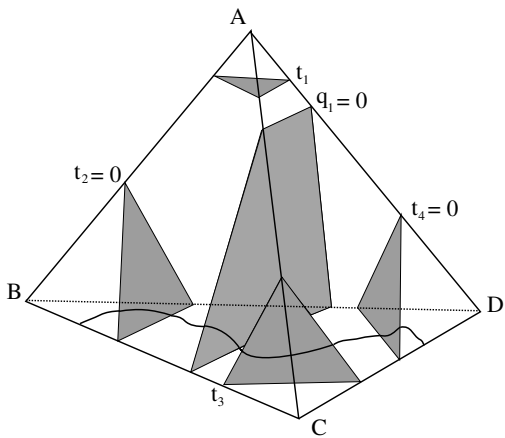}
  \caption{Modifying the normal surface equations}
  \label{modifiednormaleqns}
  \end{center}
\end{figure}

\begin{enumerate}
\item $\alpha \cap \sigma=\emptyset$:  In this case $t_2=0$, $t_3=0$, $t_4=0$, $q_1=0$, $q_2=0$, and $q_3=0$. \\
 
\item $\alpha \cap \sigma \ne \emptyset$:  Then there is a normal arc of $\alpha \cap \sigma$ intersecting $\sigma$ exactly once, say, from edge $BC$ to edge $CD$.  In this case $t_1$, $t_3$, and $q_2$ can take any value while $t_2=0$, $t_4=0$, $q_1=0$, and  $q_3=0$.  We make similar requirements for arcs that pass from edge $BC$ to $BD$ and from $CD$ to $BD$.  See Figure ~\ref{modifiednormaleqns}.\\

\end{enumerate}

%\newline
For each face $\sigma$ of the triangulation of $\partial M$ we add the above restrictions to the system of normal surface equations.  If a surface satisfies the normal surface equations as well as these additional requirements, then it is either a closed surface or its boundary consists of parallel copies of $\alpha$.  We call the new system of equations the \textit{modified normal surface equations}.  The modified normal surface equations also form a Diophantine system of linear equations.

\begin{Lem}
\label{modifieddiophantine}
The modified normal surface equations yield a finite set of  fundamental solutions such that each solution satisfying the equations can be expressed as a non-negative integer linear combination of the fundamental solutions for the modified system.  
\end{Lem}

\begin{Rem}
 In particular every incompressible Seifert surface can be expressed as a finite non-negative integer linear combination of fundamental surfaces for the modified system.
\end{Rem}

\begin{proof}
The modified normal surface equations form a Diophantine system of linear equations.  Therefore Lemma ~\ref{diophantinethm} implies that there is a finite set of non-negative integral solutions which generate using non-negative integer linear combinations, the full set of all non-negative solutions to the modified normal surface equations.  
\end{proof}

\section{Proving the main theorem}
\label{mainsect}

Throughout the remainder of this section we will assume that $K$ is a non-trivial knot in $S^3$.
Before we can prove the main theorem we need the following two lemmas.

\begin{Lem}
\label{mainlem1}
Suppose $S$ is an incompressible Seifert surface for a knot $K$ and $S=F_1+F_2+\cdots+F_k$ for some set of normal surfaces $\{F_1,...,F_k\}$ corresponding to solutions of the modified normal surface equations.  Then all of the $F_i$ are closed except for one which is a Seifert surface for the knot.  
\end{Lem}

\begin{proof}
Let $S$ be an incompressible Seifert surface for $K$ and suppose $S=F_1+F_2+\cdots+F_k$ for some set of normal surfaces $\{F_1,...,F_k\}$ corresponding to solutions of the modified normal surface equations.  Let $|\partial S|$ denote the number of components of $\partial S$.  Recall that  for each $i$, $\partial F_i$ is either empty or isotopic to copies of $\partial  S$.  For each $i$, we may isotope the components of $\partial F_i$ so that they are pairwise disjoint.  Under these conditions  $|\partial (F_1+ F_2+\cdots+F_k)| =|\partial F_1|+|\partial F_2|+\cdots+|\partial F_k|$.  Since $S$ is a Seifert surface $|\partial S|=1$.  By our hypothesis this implies $|\partial S|= |\partial (F_1+F_2+\cdots+F_k)|=|\partial(F_1)|+|\partial(F_2)|+\cdots+|\partial F_n)|=1$.  It follows that $\partial(F_i) \neq \emptyset$ for exactly one $i$, say $F_1$.  Thus $\partial S =\partial F_1$ and all of the $F_i$ are closed except for $F_1$ which is a Seifert surface.
\end{proof}

The proof of the next lemma relies on the following theorem of Jaco and Oertel \cite{JO}.

\begin{Thm}[Theorem 2.2 of \cite{JO}] 
\label{JO}
Let $T_M$ be a triangulation of the closed, irreducible 3-manifold $M$.  Suppose that $F$ is a minimal normal surface and $F=F_1+F_2$.  If $F$ is two-sided and incompressible, then both $F_1$ and $F_2$ are incompressible. 
\end{Thm}

\begin{Lem}
\label{mainlem2}
Let $K$ be a non-trivial knot in $S^3$ and $M=S^3-N(K)$.  Suppose that $S$ is an incompressible Seifert surface for $K$ such that $S=F_1+F_2$ where $F_1$ and $F_2$ satisfy the modified normal surface equations.  Then each of $F_1$ and $F_2$ are  incompressible surfaces in $M$.
\end{Lem}

\begin{proof}
By Lemma ~\ref{mainlem1} we may assume without loss of generality that $\partial{F_1}=\emptyset$ and that $F_2$ is a Seifert surface for $K$.  Let $\overline M =M\cup_T M'$ be the double of $M$, where $M'$ is a copy of $M$ and $T=\partial M$.  The double of the closed surface $F_1$ consists of two disjoint copies of itself in $\overline M$, $F_1$ in $M$, and say $F_1'$ in $M'$.  Let $\overline S$ be the double of $S$ in $\overline M$, and let $G$ be the double of $F_2$ in $\overline M$.  Then $\overline S = F_1+G+F_1'$.  Both $\overline S$ and $G$ are closed in $\overline M$.  To finish the proof of Lemma ~\ref{mainlem2} we will need the following sublemma.

\begin{Sublem}
\label{sublemforlem2}
 $\overline S$ is incompressible in  $\overline M$.
\end{Sublem}

\begin{proof}
As above, let $T=\partial M$.  Proceeding by contradiction suppose $D$ is a disk in  $\overline M$ such that $D\cap \overline S=\partial D$ and $\partial D$ does not bound a disk in $\overline S$, and $int(D) \cap S= \emptyset$.  Choose $D$ transverse to $T$ and so that $D\cap T$ has a minimal number of components.  The idea is to show that $D\cap T=\emptyset$; and therefore $D$ must be contained entirely in either $M$ or $M'$.  From here we can conclude that $D$ is a compressing disk for $S$, which is a contradiction. 

First observe that no component of $D \cap T$ is a simple closed curve.  If such a curve existed then there would be an innermost one $\alpha$ that bounds a disk $E$ in $D$ such that $E\cap T=\alpha$.  Since $T$ is incompressible there exists a disk $E'$ in $T$ such that $E^{'}\cap D=\alpha$ and $E\cup E'$ is a 2-sphere bounding a 3-ball in $\overline M$ which can be used to isotope $E'$ across $E$ removing the intersection and reducing the number of components of $D \cap T$.  Thus $D\cap T$ can have no simple closed curves.  

Now let us assume that $D\cap T$ is not empty and each component is a spanning arc of $D$.  Choose an outermost spanning arc $\alpha$ of $D \cap T$ in $D$.  Then $\alpha$ co-bounds a sub disk $E$ of $D$ to one side, and let $\beta$ be the arc in $\partial D \subset \overline S$ such that $\partial E= \alpha \cup \beta$.  The disk $E$ is contained entirely in either $M$ or $M'$.  Without loss of generality suppose $E$ is contained in $M$.  There are two possibilities.  The first is that $E$ is a boundary compressing disk for $S$.  In this case either $S$ is compressible, which is a contradiction to the fact that $S$ is an incompressible Seifert surface, or $S$ is a boundary parallel annulus, which also gives  a contradiction.  The second possibility is that $E$ is actually boundary parallel in $M$.  Here we can isotope the curve $\beta$ in $\overline S$ across the disk $E$ and past $\alpha$ in $T$ reducing the number of intersections of $D \cap T$.   So it must be the case that $D \cap T=\emptyset$.  The only possibility is that, without loss of generality, $D\subset M$,  $\partial D = D \cap S$, $int(D) \cap S=\emptyset$, and $\partial D$ is not contractible in $S$.  Thus $D$ is a compressing disk for $S$, which is a contradiction to our assumption.  This completes the proof of the sublemma.
\end{proof}

It now follows from Sublemma \ref{sublemforlem2}, and Theorem ~\ref{JO} that, since $\overline S=F_1+G+F_{1}'$ is incompressible, both $F_1$ and $G$ are incompressible.  Finally, because $G$ is incompressible in $\overline M$,  $F_2 $ is incompressible in M.  
 \end{proof}

\begin{Rem}  It follows from Lemma~\ref{mainlem1} and Lemma ~\ref{mainlem2} that one of the $F_i$ is a closed incompressible surface and the other an incompressible Seifert surface. \\
\end{Rem}

We can now prove the main theorem:
\vskip .2in

\begin{proof}[Proof of Main Theorem]

Given a non-trivial knot $K$ in $S^3$, let$M=S^3-N(K)$, and let $T_M$ be a triangulation of $M$.  Begin by fixing a preferred longitude for the knot complement, say $\alpha$.  Next, with respect to this fixed preferred longitude represented by $\alpha$, we apply Lemma \ref{curvelem} and take barycentric subdivisions of $T_M$ to obtain a new triangulation, which we will denote by $T_M'$.  This new triangulation has the property that each face of $T_M'$ on $\partial M$ meets $\alpha$ in at most one normal arc.  We are now in a position to set up the system of modified normal surface equations.

By Lemma \ref{modifieddiophantine} we can conclude that there is a finite set of fundamental solutions to the modified normal surface equations, $\{F_1,...,F_k\}$, consisting of surfaces that are closed or that have boundary isotopic to parallel copies of $\alpha$.  Any given solution to the modified normal surface equations can be expressed as a non-negative integer linear combination of these.  Now, let $\{S_1,...,S_n\}$ be the subset of the set of fundamental solutions consisting of the incompressible Seifert surface fundamental solutions, and let $\{Q_1,...,Q_m\}$ be the subset of the fundamental solutions consisting of the closed incompressible non-boundary parallel fundamental solutions.  

Given an arbitrary incompressible Seifert surface $S$, first we can isotope $\partial S$ so that it coincides with the fixed preferred longitude $\alpha$.  By Lemma \ref{fixedarclem} we can then isotope $S$ to be normal with respect to $T_M'$ while $\partial S$ remains fixed.  So by Lemma \ref{modifieddiophantine} we can write $S=a_1F_1+\cdots+a_kF_k$ where $a_1,...,a_k$ are non-negative integers.   We add the assumption that \(\sum_{i=1}^{k} a_i\) is minimal over all possible non-negative linear combination expressions for $S$.  By Lemma \ref{mainlem1} we know that all of the $F_i$ are closed except for one, which is a Seifert surface.  And by Lemma \ref{mainlem2} all of the $F_i$ are incompressible.  Note that a priori some of the surfaces in this expression for $S$ may be boundary parallel tori.  Suppose that one of the surfaces in the expression for $S$, say $F_l$, is a boundary parallel torus.  In general, taking the Haken sum of a Seifert surface with a boundary parallel torus gives a surface that is isotopic to the original Seifert surface.  That is, if $S'$ is a Seifert surface and $T$ a boundary parallel torus then $S'+T$ is isotopic to the surface $S'$.  So if one of the $F_i$ is a boundary parallel torus, then we can isotope $S$ into a position where the corresponding $a_i$ has been reduced by one.  It follows that $S=$\(\sum_{i=1}^{k} a_iF_i\) is isotopic to \(\sum_{i\neq l} a_iF_i \) $+ (a_l-1)F_l$, which is a contradiction to the minimality of \(\sum_{i=1}^{k} a_i\).  Therefore none of the $F_i$ are boundary parallel tori.  

We can conclude that all of the fundamental solutions needed in the expression for $S$ are either incompressible Sefiert surfaces and are therefore contained in the finite set $\{S_1,...,S_n\}$, or they are closed, incompressible, and non-boundary parallel, thus are contained in the finite set $\{Q_1,...,Q_m\}$.  By Lemma \ref{mainlem1} there is exactly one Seifert surface among the fundamental solutions in the expression for $S$. So we can conclude that $S$ may be isoptoed so that $S=S_i+a_1Q_1+\cdots+a_mQ_m$, where the $a_i$ are non-negative integers. 
\end{proof}

\section{AN EXAMPLE}
\label{examplesect}

In this section we present our motivating example of a knot with an infinite number of incompressible Seifert surfaces that contains a closed incompressible surface in its complement.  For a more explicit description see \cite{MSS}.  Let $K(p_1,...p_n)$ be a \emph{pretzel knot} with \emph{twist boxes} \cite{Ka} of order $p_i$.  

\begin{figure}
  \begin{center}
  \includegraphics[width=2.5in]{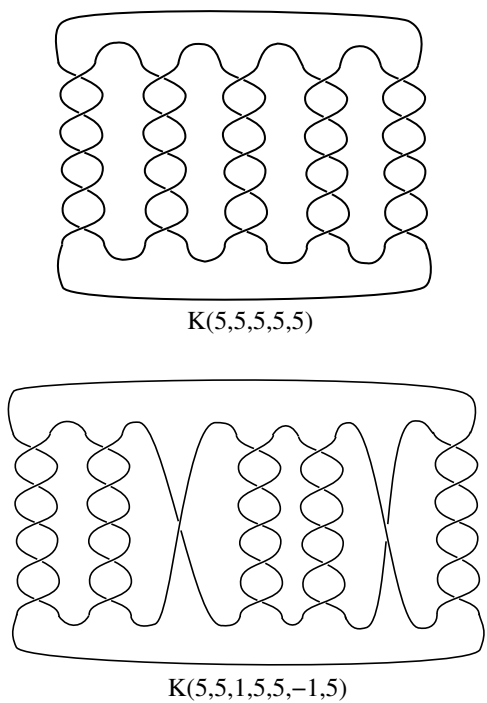}
  \caption{The pretzel knots $K(5,5,5,5,5)$ and $K(5,5,1,5,5,-1,5)$ }
  \label{pretzelflypes}
  \end{center}
\end{figure}

\begin{figure}
  \begin{center}
  \includegraphics[width=2.5in]{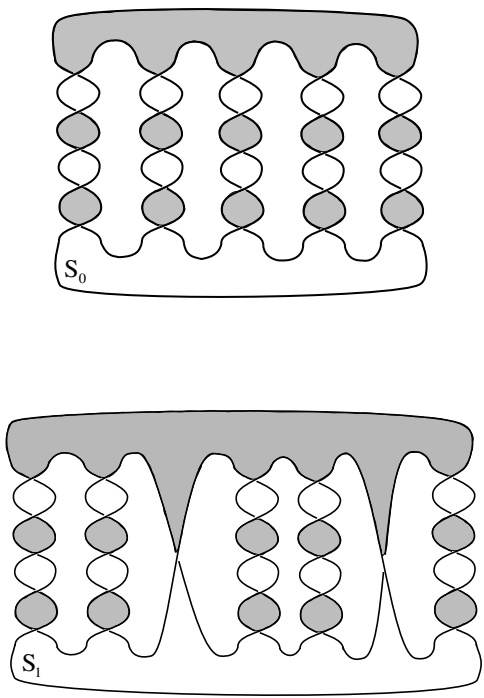}
  \caption{The Seifert surfaces $S_0$ and $S_1$}
  \label{seifert}
  \end{center}
\end{figure}

Starting with the incompressible Seifert surface $S_0$ for pretzel knot $K=K(5,5,5,5,5)$ given in Figure \ref{pretzelflypes}, one can obtain an infinite family of incompressible Seifert surfaces for $K$ by twisting the knot around a 2-sphere containing consecutive twist boxes to give a new projection for the knot $K=K(5,5,1...,1,5,5,-1,...-1,5)$.  The new projection for $K$ has a distinct new incompressible Seifert surface $S_n$, that is $S_0$ with $n$ twists.  See Figure \ref{seifert}.  Parris \cite{Pa} was the first to show that this family is incompressible. 

 Let $B$ be the 3-ball containing the third and fourth twist boxes.  Let $\Sigma=\partial B$.  Then $\Gamma={\Sigma}-N(K)$ is a four-punctured 2-sphere properly embedded in $S^3-N(K)$ with boundary consisting of meridian curves on $\partial M$.    
Let $G=\partial(B-N(K))$ be the closed surface of genus 2 constructed from $S$ by tubing together pairs of boundary components of $\partial \Gamma$.  Moriah, Sedgwick, and Schleimer \cite{MSS} show that the surfaces $S_n$ are isotopic to $S_0+nG$, and it also follows from Theorem 5.1 of \cite{MSS} that the closed surface $G$ is also incompressible.

\bibliographystyle{abbrv}

 \end{document}